\documentclass[10pt]{article}
\RequirePackage[OT1]{fontenc}
\RequirePackage{amsthm,amsmath}
\RequirePackage{natbib}
\usepackage{titlesec,epsfig,graphicx,color,amsmath,amsthm,subcaption,caption,amssymb,enumerate,multirow, setspace,hyperref}
\numberwithin{equation}{section}
\theoremstyle{plain}
\newtheorem{theorem}{Theorem}
\theoremstyle{definition}
\newtheorem{proposition}{Proposition}
\newtheorem{lemma}{Lemma}

\newtheorem{remark}{Remark}

\newcommand*\samethanks[1][\value{footnote}]{\footnotemark[#1]}

\newcommand{\round}{\lceil np \rceil}
\newcommand{\roundp}{\lceil n(1-p) \rceil}

\usepackage[left=3.6cm,right=3.6cm,]{geometry}
\usepackage{titlesec}
\titleformat{\section}{ \centering \small }{\thesection}{1 em.}{}
\titleformat{\subsection}{\small \rm \it \bf}{\thesubsection}{1 em.}{}
\titleformat{\subsubsection}{ \rm \it }{\thesubsubsection}{1 em.}{}

\begin{document}

\title{\rm {\normalsize{\bf{CAN TESTS FOR JUMPS BE VIEWED AS TESTS FOR CLUSTERS?}}}}
\author{\small{KARTHIK BHARATH}\thanks{Dept. of Statistics; The Ohio State University; 1958 Neil Avenue; Columbus, Ohio 43210.}
\and \small{VLADIMIR POZDNYAKOV}\thanks{Dept. of Statistics; University of Connecticut; 215 Glenbrook Road; Storrs, Connecticut 06269.}
\and \small {DIPAK.K. DEY}\samethanks}
\date{}
\maketitle
\begin{abstract}
We investigate the utility in employing asymptotic results related to a clustering criterion to the problem of testing for the presence of jumps in financial models. We consider the Jump Diffusion model for option pricing and demonstrate how the testing problem can be reduced to the problem of testing for the presence of clusters in the increments data. The overarching premise behind the proposed approach is in the isolation of the increments with considerably larger mean pertaining to the jumps from the ones which arise from the diffusion component. Empirical verification is provided via simulations and the test is applied to financial datasets.\\
\\
\footnotesize{KEYWORDS}: Clustering; Jump diffusions; Merton model; Test for jumps.
\end{abstract}

\section{INTRODUCTION}
 It is well-known that for an asset pricing model to circumvent arbitrage opportunities, asset prices must follow semimartingales (see \cite{DS} and \cite{HP}). In this context, jump diffusion models are popular in modeling asset prices (log asset prices, usually) in financial applications pioneered by the fundamental paper by \cite{merton}. These models are characterized by two components: a continuous component and a jump component. As described in \cite{ait6}, the continuous part of the model is present to capture the normal risk of the asset which is hedgeable whereas the jumps component can capture the default risk or news-related events. In fact, it is nowadays commonplace in applications involving high-frequency data to break up the jumps component into a large and small jump component in an effort to capture price moves which are large on the time-scale of few seconds or minutes but generally not significant on a daily scale. The focus of this article, however, will be on diffusion models with a consolidated jump component. The problem considered in this article is the following: Suppose one observes a time series of asset prices or returns over a finite length of time $[0,T]$; based on these observations, is it possible to ascertain whether the process that generated the observations comprises a jump component? The problem is of obvious importance when prediction is the primary concern. The ramifications, while constructing a model for asset pricing, of not incorporating a jump component when the underlying process which generates the data indeed does possess one, can be quite severe. The problem has received appreciable attention over the years based on several techniques; we refer to a few articles from an exhaustive list: \cite{ait1}, \cite{ait2}, \cite{carr}, \cite{barndoff1}, 
\cite{podolski} and \cite{lee}. 

~The problem can be viewed as a deconstruction problem wherein the observed series of returns are deconstructed back to their continuous and jump components. This taxonomy between the continuous and the jump components of the purported model assists us in seeking `typical' behavior of statistics based on observations under the presence and absence of the jump components. Intuitively, by constructing a test statistic which would eventually isolate the jump component under the presence of jumps in the underlying process, a suitable asymptotic hypothesis test can be developed. To elaborate, for simplicity, suppose that jumps are all positive valued and the number of jumps are finite in $[0,T]$. Based on a sampling frequency, suppose we consider the increments (difference between successive observations); we would then expect to see, primarily, two groups of data: one centered around a value which is considerably larger than the other corresponding, respectively, to the jumps and the non-jumps. Such behavior is the motivation behind the test statistics based on truncated power variations employed in \cite{ait1} and \cite{ait2}. In this article, we approach the problem of constructing a suitable test statistic through a different route: we ask if the isolation of the jumps can be viewed as a {\it model-based clustering behavior wherein the distributional properties of the model, for large samples, leads to the formation of two clusters with cluster centers far apart}. Under this setup, this article ought to be viewed as a first step towards providing a general answer applicable to a broad class of semimartingale models; while the alternative hypothesis of `jumps' encompasses a large class of models, we will focus primarily on the Merton-type model wherein the jump component is driven by a Poisson process. The Merton model (\cite{merton}) was the first model for option pricing which allowed for discontinuities in the underlying process by incorporating i.i.d. normal jump sizes with a Poisson process driving the jump process. While in recent times, very general classes of semimartingales have been employed as models for pricing, the Merton model, nonetheless, is a popular model owing to its simplicity and analytical tractability. In \cite{kou} and \cite{kou2}, the authors make a persuasive case for the Merton-type models exhorting their use in practice by using double exponential jump sizes. The double exponential jump diffusion model was calibrated and applied to market data in \cite{CT} and results were shown to be promising. Generalizations along the lines of having separate volatility and drift process albeit advantageous from a modeling perspective are not conducive for purposes of analytical tractability and interpretability. These considerations point towards use the jump diffusion model as the archetypal model in this article. 

In a recent article, motivated by the criterion function in \cite{JH2}, \cite{KB} and \cite{KB2} proposed a clustering criterion for the optimal bifurcation of a set of observations into two clusters; their approach was based on determining the point at which the data would split into clusters and this point corresponded to the zero of their criterion function. They considered sample-based versions of the clustering criterion and its zero, and proved limit theorems. We will motivate the use of clustering methods in constructing a test for jumps and consequently, employ their clustering criterion in our attempt to provide a test for jumps. In section \ref{merton} we describe the jump diffusion model and demonstrate how the testing for jumps problem can be viewed as a test for the presence of clusters. In section \ref{clustering} we detail the clustering criterion proposed in \cite{KB2} and review the relevant results from their paper. In section \ref{jumps}, we set up the hypothesis test and construct the requisite test statistic. Then, in section \ref{data} we proffer results from simulations examining the performance of the test in comparison to the test proposed by \cite{ait2} and also investigate the power of our test against particular alternatives; we then apply our test on two datasets pertaining to S\&P 500 Index returns across time windows Jan 96 - Dec 2000 and Jan 06 - Dec 10. The two time windows correspond, respectively, to periods of contrasting market behavior and it is seen that our test captures this phenomenon. Finally in section \ref{conclusion}, we summarize some of the salient features about our approach, comment on extensions and note some of its shortcomings. Proofs of results are relegated to the Appendix.  

\section{JUMP DIFFUSION MODEL}\label{merton}
We will assume that the log of the asset price follows an one dimensional It\^{o} semimartingale  process on a fixed complete probability space $(\Omega,\mathcal{F}_t,P)$, where $\{\mathcal{F}_t:t\in[0,T]\}$ is a right continuous filtration and $P$ is the data generating measure, given by
\begin{equation}\label{model}
dX_t=\mu dt+\sigma dB_t+dJ_{t_-},
\end{equation}
where the scalar $\mu \in \mathbb{R}$ represents the drift component of the process, $\sigma \in \mathbb{R}^+$, its spot volatility, $B_t$ is an $\mathcal{F}_t$ adapted standard Brownian motion and process $J_t$ is a pure jump process. For simplicity, we shall assume that $T$ is equal to 1. We will assume further that model produces observations that are collected at discrete sampling intervals $\Delta_n$ implying a regular sampling interval in $[0,1]$; we will hence suppose that we observe $X_t$  at $n$ discrete times $0 \leq \Delta_n \leq 2\Delta_n \leq \cdots \leq n\Delta_n \leq 1$ where $\Delta_n=1/n$. Our intention is to study the model as $\Delta_n \rightarrow 0$ as $n\rightarrow \infty$. Based on a discretely observed trajectory or path of $X$, our objective is to assign the observed path to two complementary sets:
\begin{align}\label{sets}
\Omega^j&=\bigg\{\omega: t \to X_t(\omega) \text{ contains jumps in }[0,1]\bigg\},\nonumber \\
\Omega^c&=\bigg\{\omega: t \to X_t(\omega) \text{ does not contain jumps in }[0,1]\bigg\}.
\end{align}
If we choose $\Omega^c$, we are implicitly stating that we are choosing $X_t=X_0+\mu t+\sigma B_t$ with a.s. continuous paths. For technical reasons and also as a natural way to assess if the process $X$ has jumped in $[0,1]$, it is common to consider the increments
\begin{equation}\label{increments}
W_i=X_{(i+1)\Delta_n} - X_{i\Delta_n} \quad 0 \leq i \leq n-1.
\end{equation}
If the jump process $J_t$  is assumed to be a Levy process, then by the independent increments property of $B_t$ and $J_t$, 
$W_1,\ldots,W_{n-1}$ are i.i.d. random variables. The trick usually used is to construct a statistic based on $W_i$, independent of the model parameters, $\mu, \sigma$ and any additional parameters for $J_t$,  in such a way that, as $n \to \infty$, the statistic's behavior would be dominated by the jumps component if $X_t$ does indeed jump in $[0,1]$. For instance, in \cite{ait6}, realized power variations of $W_i$ suitably truncated, of the form
\[
U(p,u_n,\Delta_n)=\displaystyle \sum_{i=1}^{n-1}|W_i|^p \mathbb{I}_{|W_i| > u_n},
\]
where $u_n$ is a deterministic sequence of truncation levels which tend to 0 as $n \to \infty$, were considered; typically $u_n$ is taken to be a function of $\Delta_n$ with some other constants not depending on $n$. The idea is that by a judicious choice of $p$, one can eliminate all the increments which correspond to the continuous part of the model to end up with a value for $U(p,u_n,\Delta_n)$ completely dominated by the jump component. The power $p$ plays an important role in the behavior of $U(p,u_n,\Delta_n)$: for $p>2$, the contributions from the jump component dominates the value of $U(p,u_n,\Delta_n)$ whereas for $p<2$, the continuous part dominates. Therefore, a finite value for $U(p,u_n,\Delta_n)$ for a certain choice of $p$, or some functional thereof, can be used to construct a test for jumps. This seemingly simple but powerful idea is exploited in a different setting in this article. 

In contrast to the approach adopted by \cite{ait6}, it can be intuited that such a forced separation of the increments can be performed by suitably choosing a truncation level too. {\it Our approach is characterized, in a certain sense, by the determination of truncation level $u_n$, which splits the observations into two clusters: one pertaining to the continuous part and the other to the jump part}. Our test statistic will be based on the truncation level at which the two groups clearly bifurcate. To elucidate, suppose we assume that the process $J_t$ is a finite-activity jump process like a Compound Poisson. To aid intuition, suppose additionally that the jump-sizes are all positive. Let us consider the complementary statistics:
\[
U_1(u_n,\Delta_n)=\frac{1}{(n-1)F_n(u_n)}\displaystyle \sum_{i=1}^{n-1}W_i \mathbb{I}_{W_i \leq u_n};  U_2(u_n,\Delta_n)=\frac{1}{(n-1)(1-F_n(u_n))}\displaystyle \sum_{i=1}^{n-1}W_i \mathbb{I}_{W_i > u_n},
\] 
where $F_n$ is the empirical distribution function of the $W_i$. The two statistics correspond to the average of the observations greater and lesser than the truncation level. For an `optimal' choice of a truncation level $u_n$, one would be able to observe a clear separation between the large values of $W_i$ (jumps) and the smaller values (continuous part); the optimal level would have to be determined by comparing \emph{all} possible averages of the two groups. In view of this, we ask the question `Can the optimal truncation level be determined for the model with no jumps. i.e. the model comprising just the Brownian motion with drift?'. An answer in the affirmative would point towards the usage of an estimate of the optimal truncation level as a candidate for a test statistic for the test for jumps. 

The preceding discussion, from a statistician's perspective at least, would suggest to the use of a mixture model as suitable probabilistic mechanism as an explanation for observed increments $W_1,\ldots, W_{n-1}$---in particular, we are interested in a two component mixture model with a component each for the continuous and the jump parts. Clustering methods have been widely used to analyze mixture models and perform statistical inference. As a first step towards developing a test for a large class of semimartingale models with a Brownian component, in this article we consider the Merton model or its close cousins and demonstrate how clustering methods can be used in constructing a test statistic for jumps. We believe that a good understanding of the simple yet fairly general and popular Merton model would be a step in the right direction towards developing tests for jumps using clustering criteria. 

\subsection{Merton model}
\begin{figure}[!h]
\begin{center}
\includegraphics[scale=0.4]{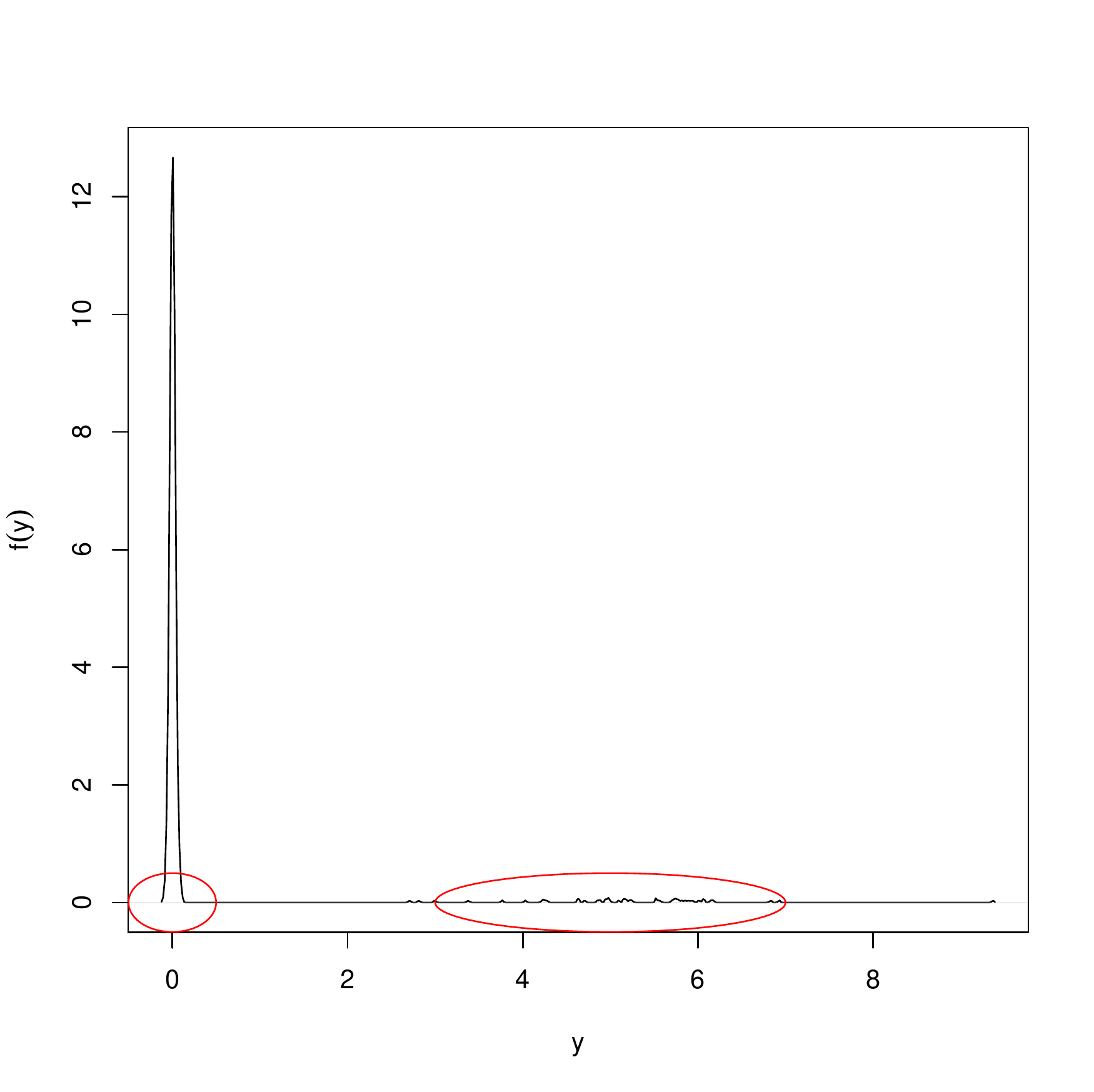}
\caption{\small{The figure shows the density $f$ for the Merton model observed discretely at $n=1000$ points with $\mu=2$, $\sigma^2=1$, $\lambda=5$ and jump sizes taken as i.i.d. normal random variables with mean 5 and variance 1. The clustering of the points around 0 and around $\lambda=5$ is quite clear and is denoted by the red ellipses.}}
\label{split_fig}
\end{center}
\end{figure}

Under the setup of the model in (\ref{model}), the Merton model for assets pricing $X_t$, is given by,
\begin{equation}\label{merton_model}
X_t=X_0+\mu t+\sigma B_t+\displaystyle \sum_{k=0}^{P_t(\lambda)}J_k \quad 0\leq t \leq T,
\end{equation}
where $P_t$ is a Poisson jump process with intensity $\lambda$ with jumps sizes represented by i.i.d. random variables $J_k$. It is assumed that $B_t$, $P_t(\lambda)$ and $\{J_k\}$ are independent and all parameters are assumed to be unknown.  
For the purposes of demonstrating the validity our method, we consider a special case of (\ref{merton_model}) when all the jumps are of unknown constant size $h>0$. While this condition is perhaps not realistic in applications, as a first step towards understanding the utility of clustering to the testing problem, it is reasonable. As a justification, consider the two sets in (\ref{sets}). Indeed, if our test rejects the null hypothesis of no jumps, we assign the observed path of $X$ to $\Omega^j$, a very general class and is independent of our assumption on the nature of the jumps. Obviously, the general nature of the alternative hypothesis might not guarantee good power against a very specific subset of processes with jumps. Nevertheless, if the primary motivation is to just test for the presence or absence of jumps based on discrete observations, the level of any proposed test is independent of the assumed structure of the nature of the jump component. The test proposed, therefore, in a certain sense, would be an omnibus test. The circumscription of our approach in this article to the case of constant positive size jumps is hence reasonable bearing in mind the preceding discussion. Additionally, when jumps are of size $h>0$, the model in (\ref{merton_model}) reduces to 
\begin{equation}\label{our_model}
X_t=X_0+\mu t+\sigma B_t+hP_t(\lambda) \quad 0\leq t \leq T.
\end{equation}
We consider the increments
\begin{equation}\label{increments}
W_i=X_{(i+1)\Delta_n} - X_{i\Delta_n} \quad 0 \leq i \leq n-1,
\end{equation}
the discretized version of $X$, which, as a consequence of the independent increments property of the Brownian motion and the Poisson process, are i.i.d. with density
\begin{equation}\label{original_density}
f(w) = \displaystyle \sum_{k=0}^\infty \phi_{\mu^*+hk,\sigma^*}(w)\frac{e^{-\lambda^*}\lambda^{*k}}{k!},
\end{equation}
where $\mu^*=\mu \Delta_n$, $\sigma^*=\sigma \sqrt{\Delta_n}$, $\lambda^*=\lambda \Delta_n$ and $\phi_{a,b}$ represents the density of a normal random variable with mean $a$ and standard deviation $b$. The density is an infinite mixture of Gaussian distributions with the mixing proportions obtained from a Poisson distribution with rate $\lambda^*$. We thus have
\begin{equation*}
f(w)= \phi_{\mu^*,\sigma^*}(w)e^{-\lambda^*}+\phi_{\mu^*+h,\sigma^*}(w)e^{-\lambda^*}\lambda^*+\displaystyle \sum_{k=2}^\infty \phi_{\mu^*+hk,\sigma^*}(w)\frac{e^{-\lambda^*}\lambda^{*k}}{k!}.
\end{equation*}
Since $\Delta_n\rightarrow 0$ as $n \rightarrow \infty$ and from the orderliness of the Poisson process, $P(P_{t+\Delta_n}-P_t>1|P_{t+\Delta_n}-P_t \geq 1)$ goes to $0$ as $n \rightarrow \infty$; the probability of observing two or more jumps in one increment goes to zero as $n\rightarrow \infty$. This implies that we can reduce the infinite mixture density to a two component mixture density since for large $n$, $e^{-\lambda^*} \asymp (1-\lambda^*)$ and $e^{-\lambda^*}\lambda^*\asymp \lambda^*$, we have
 \begin{equation*}
 f(w)= (1-\lambda^*)\phi_{\mu^*,\sigma^*}(w)+\lambda^*\phi_{\mu^*+h,\sigma^*}(w)+o(n^{-1}).
 \end{equation*}
Consequently, when $n$ is large, for our purposes, we can assume that the increments $W_i, 1\leq i \leq n-1$, are i.i.d. from density
 \begin{equation}\label{density}
 g(w)= (1-\lambda^*)\phi_{\mu^*,\sigma^*}(w)+\lambda^*\phi_{\mu^*+h,\sigma^*}(w).
 \end{equation}
For large $n$, since $\mu^*$, $\sigma^*$ and $\lambda^*$ are very close to zero, we would expect to see most of the observations concentrated around 0 and few clustered around $h$. The unknown $h$, in a certain sense, is related to the intensity $\lambda$ of the Poisson process typified by the observation that as $n$ increases the jumps concentrate around their expected value. Therefore, for large sample sizes, intuitively we would expect the density $g$ in (\ref{density}) to almost have point mass at zero and very little mass $h$ or $\lambda$; Figure \ref{split_fig} illustrates this behavior. 
Starting from observations from density from (\ref{original_density}) we have made the transition to observations from (\ref{density}). The relevant density is a two-component Gaussian mixture wherein the component variances and the mixing proportions are going down to zero and one as sample size increases. The problem of testing for jumps based on a discretely observed process can hence be reduced to a statistical problem of testing if the given data is indeed a random sample obtained from a mixture of two normal populations wherein one of the mixing proportions and component variances tend to zero as sample size approaches infinity. This type of mixture is different from classical mixture-models since the weights in our model tend to zero and one with increasing sample size; existing clustering methods for testing in mixture-models are rendered inapplicable in our setup. 

\section{A SUITABLE CLUSTERING CRITERION}\label{clustering}
 Bearing in mind the density in (\ref{density}), the objective is to construct a test statistic by estimating the optimal truncation level which would bifurcate the observations corresponding to the continuous and the jump parts; furthermore, the truncation level ought to be independent of the parameters $(\mu,\sigma^2, \lambda,h)$ since we do not wish to estimate them thereby compromising on the power. As a consequence, we consider the nonparametric clustering criterion proposed in \cite{KB2}. Their criterion is based on determining the point at which data is broken up into clusters---this is precisely what is required for our test for jumps. In this section, we review some preliminaries of the clustering criterion proposed with the view of constructing a test statistic that can be used to test for the presence of clusters in a set of observations; in other words, we would construct a test for the presence of jumps.

While it is tempting to recast the general setup in \cite{KB2} in terms of our specific problem of data from a normal mixture, it is important to understand why their criterion function can be useful at all in determining the truncation level. We thus adopt the general setup used in their article and comment, wherever necessary, on the adaptation of their results to the testing problem. It is to be noted, however, that their results in their existing form are not amenable for direct application in our setting. This point will be elucidated in the subsequent sections; we are hence tasked with suitably modifying their results to tailor our needs.

Suppose $W_1,W_2,\cdots,W_n$ are i.i.d. random variables with cumulative distribution function $F$. Denote by $Q$ the quantile function associated with  $F$. We make the following assumptions:
\begin{description}
\item [$A1$.] $Q$ is the unique inverse for $0<p<1$ and $F$ is absolutely continuous with respect to the Lebesgue measure with density $f$.
\item [$A2$.] $E(W_1)=0$ and $E(W_1^2)=1$.
\item [$A3$.] $Q$ is twice continuously differentiable at any $0<p<1$.
\end{description}
\begin{remark}
Note that assumptions $A1$ and $A3$ are satisfied by the distribution function of a normal random variable; assumption $A2$ will be clarified soon. This is relevant since the test for presence of jumps has been reduced to the test for the presence of clusters in a sample from a normal population(s).
\end{remark}
\subsection{Empirical Cross-over Function}
The \emph{cross-over} function, for $0<p<1$, is defined as
\begin{equation}\label{crossover}
G(p)= \frac{1}{p}EW_1\mathbb{I}_{W_1\leq Q(p)}+\frac{1}{1-p}EW_1\mathbb{I}_{W_1> Q(p)}-2Q(p).
\end{equation}
A point $p_0$ which solves $G(p)=0$ is referred to as the \emph{split point}. 
\begin{remark}
The cross-over function $G$ is a function of the derivative of 
\[
B(Q,p)= \frac{1}{p}\left[EW_1\mathbb{I}_{W_1\leq Q(p)}\right]^2+\frac{1}{1-p}\left[EW_1\mathbb{I}_{W_1> Q(p)}\right]^2-E^2W_1,
\]
referred to as the \emph{split function} in \cite{JH2}. The function $B$ can be viewed as the between cluster sums of squares and the point $p_0$ at which $B(Q,p)$ is maximized (with respect to $p$) coincides with the zero $G$. \cite{JH2} considered sample versions of $B$ and $p_0$ and investigated their asymptotic behavior. A test statistic for test for bimodality was constructed using the sample version of $p_0$. 
\end{remark}
Let us briefly examine the behavior of the cross-over function $G$: It starts positive, crosses zero and assumes negative values; the point of crossing is of chief interest and its location in $(0,1)$ can be used as an indication of the nature of the underlying distribution: symmetric, unimodal or not unimodal. 
\begin{figure}
\centering
\begin{subfigure}[ht]{0.45\textwidth}
\centering
\includegraphics[width=\textwidth]{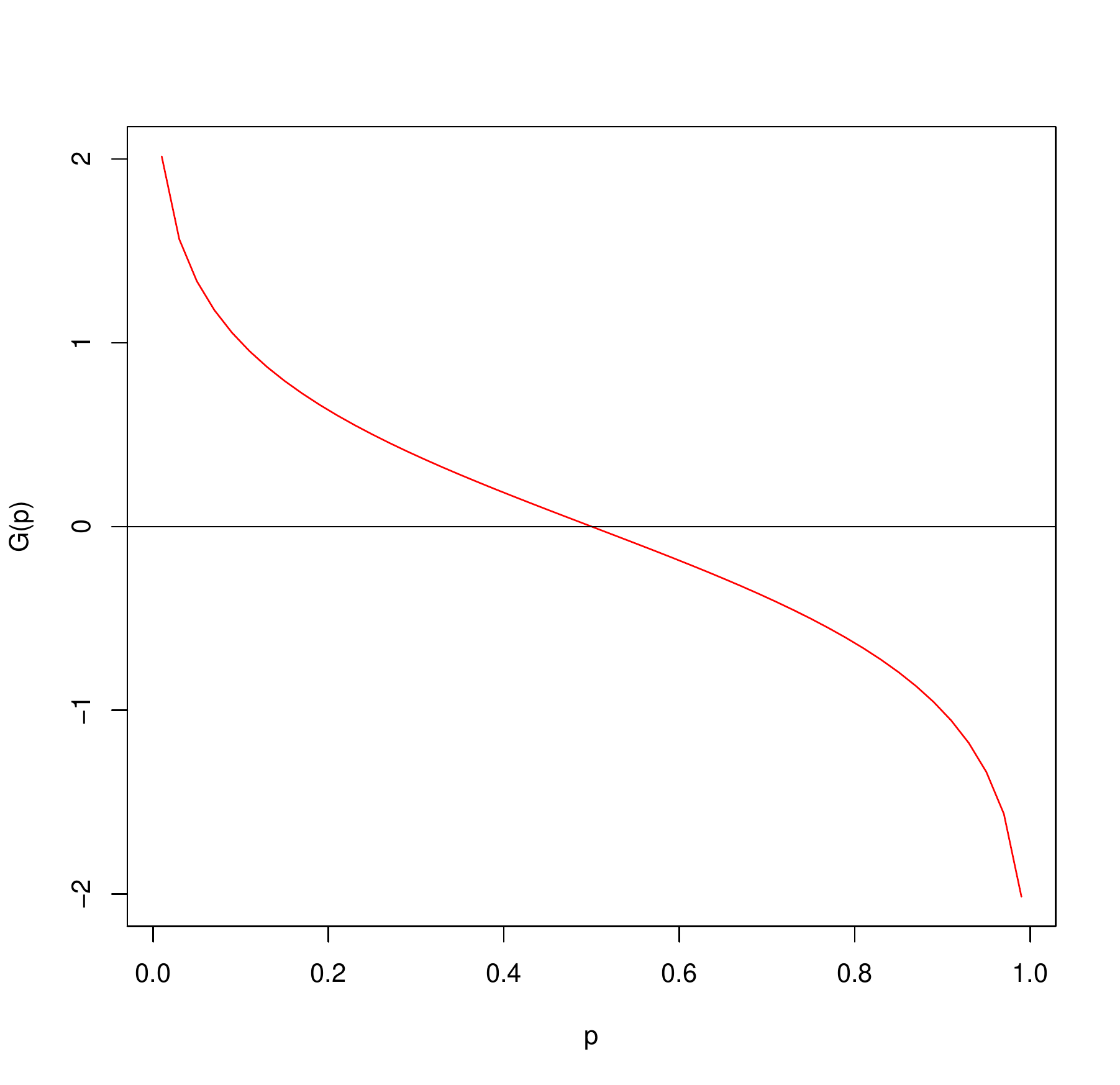}
\caption{\label{fig2a}\small{$G$ for $N(0,1)$, $p_0=0.5$.}}

\end{subfigure}
\centering
\begin{subfigure}[ht]{0.45\textwidth}
\centering
\includegraphics[width=\textwidth]{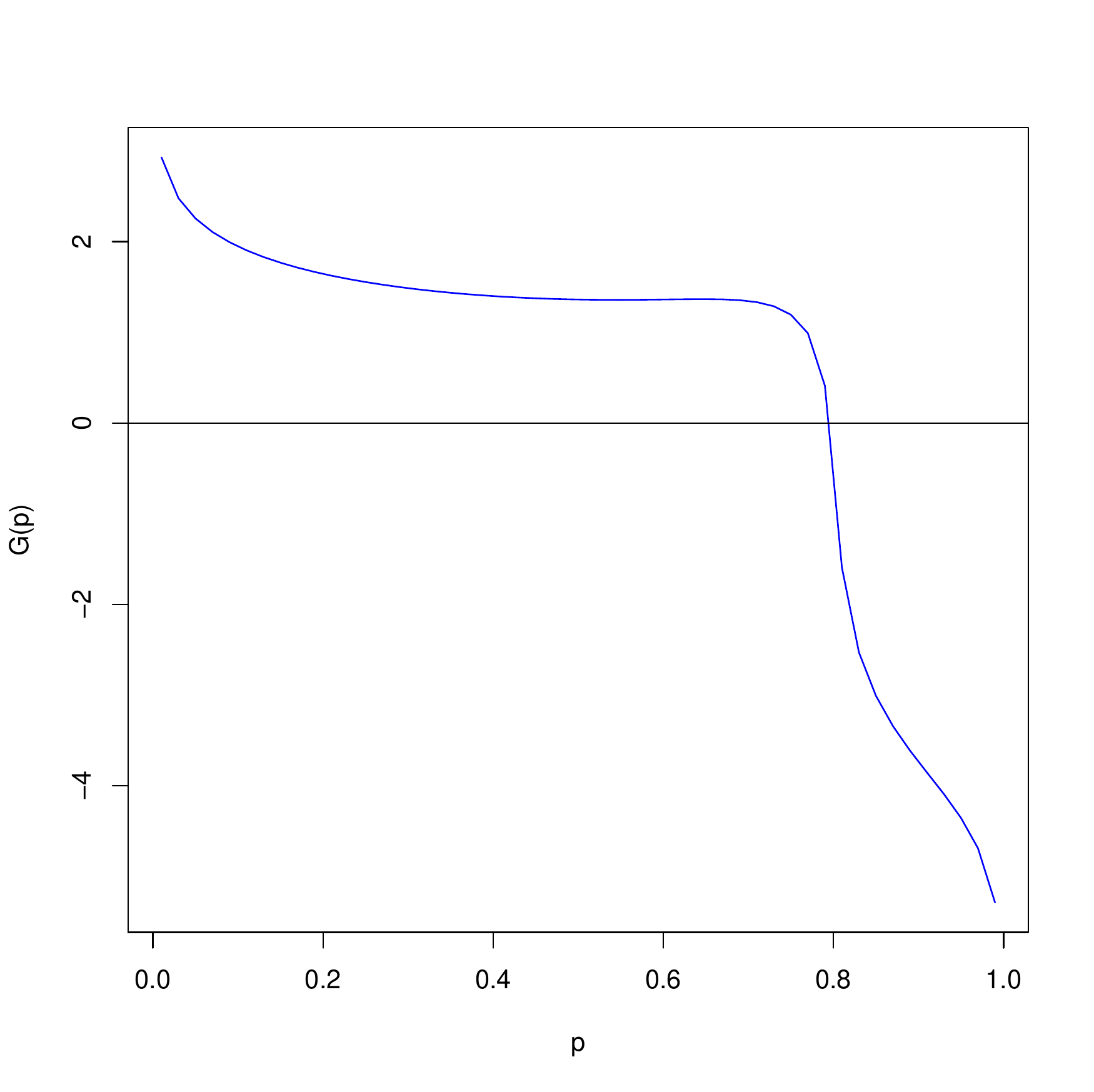}
\caption{\label{fig2b} \small{$G$ for $0.8 N(0,1)+0.2 N(5,1)$, $p_0\approx 0.8$.}}

\end{subfigure}
\end{figure}
In particular, it is easy to see that for the standard normal density $0.5$ is the split point (Figure \ref{fig2a}). What is important is that in the case of a {\it mixture of two normals with same variance but differing means, the split point is not 0.5}. This can be clearly inferred from the $0.8 N(0,1)+0.2 N(5,1)$ density in Figure \ref{fig2b} where the split point $p_0$ is approximately 0.8. These considerations suggest using the split point as an indicator for the presence or absence of clusters. One may question the behavior of the $G$ and the location of the split point for an entire parametric class of normals; fortunately, by virtue of the definition of $G$, the split point $p_0$ is invariant to scaling and translations of the standard normal density. This fact is very convenient in our setup of testing for jumps as we can then disregard the estimation of the parameters $\mu$, $\sigma^2$ and $h$ in our attempt to construct a suitable test statistic. 

Since the cross-over function $G$ appears to be a good candidate for use in the testing problem, we now review its empirical counterpart proposed in \cite{KB2}. Suppose $W_{(1)},\ldots, W_{(n)}$ are the order statistics corresponding to the i.i.d. observations $W_i,1 \leq i \leq n$. The \emph{Empirical Cross-over Function} (ECF) is then defined as
\begin{equation}\label{G_function}
 G_n(p)=\frac{1}{k}\sum_{j=1}^k W_{(j)}-W_{(k)} +\frac{1}{n-k}\sum_{j=k+1}^n W_{(j)}-W_{(k+1)},
\end{equation}
for $\frac{k-1}{n}\leq p <\frac{k}{n}$
{ and
\begin{equation}\label{G_function*at*1}
 G_n(p)=\frac{1}{n}\sum_{j=1}^n W_{(j)}-W_{(n)},
\end{equation}
for $\frac{n-1}{n}\leq p <1$,  where $1 \leq k \leq n-1$.}
\begin{remark}
Observe that for a fixed $0<p<1$
\begin{align*}
\frac{1}{k}\sum_{j=1}^k W_{(j)}-W_{(k)} &= \frac{1}{\round}\sum_{j=1}^{\round} W_{(j)}-W_{(\round)}, \\
\frac{1}{n-k}\sum_{j=k+1}^n W_{(j)}-W_{(k+1)} &= \frac{1}{\roundp}\sum_{j=\round+1}^{n} W_{(j)}-W_{(\round +1)},
\end{align*}
where $\lceil x\rceil$ represents the smallest integer not less than $x$. For a fixed $p \in (0,1)$, the sums shown above are \emph{trimmed} sums. More precisely, since $\frac{\round}{n} \rightarrow p$ and $\frac{\roundp}{n} \rightarrow 1-p$, they represent the case of \emph{heavy} trimming. We will employ the two notations interchangeably when there is no confusion. 
\end{remark}
The random quantity $G_n$, represents the empirical version of the cross-over function $G$ and determines the split point for the given data. The intuition behind the sample criterion function is simple: $G_n$ is based on the distances between the sample quantiles $W_{(k)}$ and $W_{(k+1)}$ and the means of the observations lesser and greater than them respectively. Then, by checking the distances for all possible indices $k$, one hopes to ascertain the particular $k^*$ at which the distances match up, viz., the function $G_n$ becomes $0$. Based on that particular index $k^*$, one is then able to infer if the sample perhaps was obtained from a population with `more than one mean', and estimate the split point. Notice the striking similarity between the rationale employed here and the one proposed while considering the statistics $U_i(u_n,\Delta_n)$ for $i=1,2$ for determining the optimal truncation level $u_n$. Indeed, the statistics $U_i(u_n,\Delta_n)$ are based on truncated sums as opposed to the trimmed sums used in the definition of $G_n$. This difference is not of great significance; \cite{KB3} employed an alternative definition of $G_n$ using truncated sums and obtained very similar asymptotic results. 

It is noted in \cite{KB} that the ECF $G_n$ is an L-statistic with irregular weights and hence not amenable for direct application of existing asymptotic results for L-statistics.
Observe that
\begin{align*}
 G_n\left(\frac{0}{n}\right) &= W_{(1)}-W_{(1)}+\frac{1}{n-1}\sum_{j=2}^n W_{(j)}-W_{(2)}\quad \geq 0,\\
  G_n\left(\frac{n-1}{n}\right)&= \frac{1}{n}\sum_{j=1}^{n} W_{(j)}-W_{(n)} \quad \leq 0.
\end{align*}
This simple observation captures the typical behavior of the empirical cross-over function. As is the case with $G$, it starts positive and then at some point crosses the zero line.  The index $k$ at which this change occurs determines the datum $W_{(k)}$ at which the split occurs. In \cite{KB}, it is shown that $G_n(p)$ is a consistent estimator of $G(p)$ for each $0<p<1$ and a functional CLT was also proved for $\sqrt{n}(G_n(p)-G(p))$ for $p \in [a,b]$ with $0<a<b<1$; it was shown that the limit stochastic process was Gaussian.
\subsection{Empirical split point}
Our interest, however, is in the split point $p_0$. On the scale of the random variables, the optimal truncation level is $Q(p_0)$. Since this quantity is unknown, we would like to estimate it using the data; this leads us to the \emph{empirical split point} defined in \cite{KB2}. There is a technical issue here: in \cite{KB2} the empirical split point is defined on $[a,b]$ where $0<a <b<1$ and all the asymptotic results were proved on $[a,b]$. In our application, the jumps or the clusters, for large $n$, are close to the boundary near 1. In order to use the results in \cite{KB2}, we would first need to amend their definition of the empirical split point to $(0,1)$ when the underlying distribution is normal and then ensure that their results are still valid in our setting. We shall demonstrate that this is possible and we therefore redefine the empirical split point as 
\begin{equation*}
 p_n:=\left\{
 \begin{array}{l}
 0, \mbox{ if } G_n\left(\frac{k-1}{n}\right) < 0 \quad \forall 1 \leq k \leq n-1;\\
 \\
 1, \mbox{ if } G_n\left(\frac{k-1}{n}\right) > 0 \quad \forall 1 \leq k \leq n-1;\\
 \\
 \frac{1}{n}\left[\max\{1 \leq k \leq n-1: G_n\left(\frac{k-1}{n}\right)G_n\left(\frac{k}{n}\right) \leq 0\}\right], \mbox{ otherwise.}
 \end{array}
 \right.
\end{equation*}
Notice that we are defining the empirical split point to be the index $k$ immediately after which $G_n$ becomes negative. Implicitly, the assumption here is the $G_n$ upon crossing the line $y=0$, can only perhaps crossover again within a small neighborhood of true split point $p_0$; such as assumption is required while proving asymptotic properties of $p_n$. In general, such an assumption is difficult to verify and conditions on the distribution function which guarantee such a behavior is not clear; this was noted in \cite{KB2} and they hence restricted their results to $[a,b]$ with $0<a<b<1$. We extend their results for $0<p<1$ using the following lemma:
\begin{lemma}\label{supGn}
Suppose $F$ is the standard normal distribution function. Then, there exists a $0<b<1$ such that, as $n \to \infty$,
\[
\sup_{p>b}G_n(p) \overset{P}\to -\infty.
\]
\end{lemma}
\begin{remark}
Lemma \ref{supGn} is imperative for the employment of $G_n$ and $p_n$ in the problem of testing for jumps since the cluster corresponding to the jumps can, in principle, be formed very close to $p=1$; we do not have the luxury of restricting ourselves to a closed sub-interval of $(0,1)$. Since the ECF $G_n$ tends to $-\infty$ as  $p \to 1$ we are assured that $G_n$ would do crossover the line $y=0$ far away from the true split $p_0$. Thus $G_n$, in conjunction with $p_n$, can be used as a tool to develop a test for the presence of clusters by examining the asymptotic behavior of $p_n$ under the null hypothesis of no-jumps. Indeed, the fact that our sample is from a normal distribution assists us in the proof. However, it is pertinent to note that the result of Lemma \ref{supGn} is applicable under a broader setup for distribution functions which have tail behavior similar to the normal. 
\end{remark}
\begin{remark}
From a practical perspective, what is important is that $p_n$ is invariant to scaling and translations of the data. Notice that if for constants $\alpha>0$ and $\beta$ and $i=1,\dots,n$,
\begin{equation*}
Z_i=\alpha W_i+\beta  ,
\end{equation*}
and we define $G_n^z$ to be the ECF based on $Z_i$, then,
\begin{align*}
G_n^z\left(\frac{k-1}{n}\right)    &= \frac{1}{k}\displaystyle \sum_{j=1}^k Z_{(j)}-Z_{(k)} + \frac{1}{n-k}\displaystyle \sum_{j=k+1}^n Z_{(j)}-Z_{(k+1)}\\
        &= \alpha \left[\frac{1}{k}\displaystyle \sum_{j=1}^k W_{(j)}-W_{(k)} + \frac{1}{n-k}\displaystyle \sum_{j=k+1}^n W_{(j)}-W_{(k+1)}\right]\\
        &= \alpha G_n\left(\frac{k-1}{n}\right),
\end{align*}
and therefore, $G_n^z$ and $G_n$ cross-over $0$ at the same point; this shows that $p_n$ is invariant to scaling and translations. This is of primary importance to us while constructing the test for jumps since it frees us from having to estimate the drift and the volatility coefficients and provides the rationale behind assumption $A2$.
\end{remark}
How accurate is the estimate $p_n$ of the split point $p_0$? The following two theorems, which are similar to their counterparts in \cite{KB2}, shed light on this issue. We essentially extend their results for $p$ in $[a,b]$ with $0<a<b<1$ to all $p \in (0,1)$. The proofs for the theorems carry over with minimal change from theirs assisted by Lemma \ref{supGn}. We hence just state them under the our modified setup and omit the proofs.  
\begin{theorem}
\label{th1}
Assume $A1-A3$ hold. Suppose that $G(p)=0$ has a unique solution, $p_0$. Then for any $0<p_0<1$
\[
p_n\overset{P}\rightarrow p_0,
\]
as $n\rightarrow \infty$.
\end{theorem}
The empirical split point is shown to be a consistent estimator of the theoretical split point; however, for constructing a test, we need more. We require the nature of the deviation of $p_n$ from $p_0$ and a Central Limit Theorem (CLT) is proved in \cite{KB2}. Before stating their theorem, define

\begin{align*}
\theta_p=&\phantom{+}\frac{1}{p}W_1\mathbb{I}_{W_1<Q(p)}-\frac{1}{p}Q(p)\mathbb{I}_{W_1<Q(p)}\\
         &+\frac{1}{1-p}W_1\mathbb{I}_{W_1\geq Q(p)}-\frac{1}{1-p}Q(p)\mathbb{I}_{W_1\geq Q(p)}\\
         &+\frac{2\mathbb{I}_{W_1<Q(p)}}{f(Q(p))}.
\end{align*}
Note that, 
\begin{align}\label{gprime}
G^{'}(p)&=\frac{1}{p}\left[Q(p)-\frac{EW_1\mathbb{I}_{W_1 \leq Q(p)}}{p}\right]
-\frac{1}{1-p}\left[Q(p)-\frac{EW_1\mathbb{I}_{W_1>Q(p)}}{1-p}\right]-2Q'(p).
\end{align}
\begin{theorem}\label{normal}
Assume $A1-A3$ hold. Suppose that $G(p)=0$ has a unique solution, $p_0$, and $G'(p_0)<0$.  Then, as $n\rightarrow \infty$,
\[
\sqrt{n}(p_n-p_0)\Rightarrow N\left(0,\frac{Var(\theta_{p_0})}{G'^2(p_0)}\right).
\]
\end{theorem}
Theorem \ref{normal}, theoretically,  provides us with a test statistic suitable to for the presence of clusters or jumps. However, the asymptotic variance involves population quantities $Var(\theta_{p_0})$ and $G'(p_0)$ which are unknown. We will provide a consistent estimator for the asymptotic variance which can then be used to develop the test. 
\begin{proposition}\label{variance}
Let 
\begin{align*}
\eta_n(p_n)&=\frac{S_{nl}}{p_n}+\frac{W^2_{(\lceil np_n\rceil)}}{p_n}+\frac{S_{nu}}{1-p_n}+\frac{W^2_{(\lceil np_n\rceil)}}{1-p_n}+4p_n\hat{Q}^{'}(p_n)\\
&\quad - \frac{2W_{(\lceil np_n\rceil )}T_{nl}}{p_n}-\frac{2W_{(\lceil np_n\rceil)}T_{nu}}{1-p_n}+4T_{nl}\hat{Q}^{'}(p_n)\\
&\quad -4W_{(\lceil np_n\rceil )}\hat{Q}^{'}(p_n)
 \quad -\left[T_{nl}+T_{ul}-2W_{(\lceil np_n\rceil)}+2p_n\hat{Q}^{'}(p_n)\right]^2,
\end{align*}
where
\begin{align*}
S_{nl}=\frac{\displaystyle \sum_{i=1}^{\lceil np_n\rceil}W^2_{(i)}}{\lceil np_n \rceil},\quad S_{nu}=\frac{\displaystyle \sum_{i=\lceil np_n\rceil +1}^{n}W^2_{(i)}}{\lceil n(1-p_n) \rceil},
\quad T_{nl}=\frac{\displaystyle \sum_{i=1}^{\lceil np_n\rceil}W_{(i)}}{\lceil np_n \rceil},\quad T_{nu}=\frac{\displaystyle \sum_{i=\lceil np_n\rceil +1}^{n}W_{(i)}}{\lceil n(1-p_n) \rceil},
\end{align*}
and 
$$\hat{Q}^{'}(p_n)=n\big(W_{(\lceil np_n\rceil +1)}-W_{(\lceil np_n \rceil)}\big).$$
Also, let 
\begin{align*}
\delta_n(p_n)=\frac{1}{p_n}\left[W_{(\lceil np_n \rceil)}-\frac{1}{p_n}T_{nl}\right]-\frac{1}{1-p_n}\left[W_{(\lceil np_n \rceil)}-\frac{1}{1-p_n}T_{nu}\right]-2\hat{Q}^{'}(p_n).
\end{align*}
Then, as $n \to \infty$,  
\[
\frac{\eta_n(p_n)}{\delta^2_n(p_n)} \overset{P} \to \frac{Var(\theta_{p_0})}{G'^2(p_0)}.
\]
\end{proposition}
Proposition \ref{variance} can be used in conjunction with Theorem \ref{normal} to construct a test for clusters via Slutsky's theorem.
\section{TEST FOR JUMPS}\label{jumps}
For the i.i.d. increments $W_1,\ldots, W_{n-1}$ from density (\ref{density}), it appears that the order statistic $W_{(\lceil np_n\rceil)}$ corresponds to the estimate of the optimal truncation level which would separate the increments into the continuous and the jumps clusters. However, it is not immediately clear from the asymptotic results for clustering, why $W_{(\lceil np_n\rceil )}$ would be a good choice bearing in mind the rather unusual nature of the density in (\ref{density}): a mixture distribution with some of the parameters tending to zero with increasing sample size. The key question to be addressed is the following: why is it the case that when $n$ is large and there is a clear separation between two clusters---one corresponding to the Brownian component around 0 and another around $h$--- the Empirical Cross-over Function (ECF) captures it? We will prove a theorem showing why the ECF is useful in our setting; Lemma \ref{supGn} is key in this regard and ensures that the jumps, which are observed close to $p=1$, can be captured accurately. We first provide an informal explanation as to why the theorem is reasonable.  

Recall that for large $n$, we can assume that the density of the i.i.d. (to be accurate, they form a triangular array) increments $W_i, i \leq i \leq n-1$ is  
\begin{equation*}
 g(w)= (1-\lambda^*)\phi_{\mu^*,\sigma^*}(w)+\lambda^*\phi_{\mu^*+h,\sigma^*}(w) \quad w \in \mathbb{R},
 \end{equation*}
where $\mu^*=\mu \Delta_n$, $\sigma^*=\sigma \sqrt{\Delta_n}$, $\lambda^*=\lambda \Delta_n$ and $\Delta_n=1/n$.
Using the results from Section \ref{clustering}, suppose we were to use the ECF and determine the empirical split point $p_n$. Let $k^*=\lceil np_n\rceil$, the index after which the ECF turns negative. Note that for the density $g$, for large $n$, there is a separation of approximately $h$ between the adjacent order statistics, $W_{(k^*)}$ and $W_{(k^*+1)}$. Since the variances of the clusters are tending to zero together, intuition tells us that $(W_{(k^*)}-W_{(k^*+1)})$ should be approximately $-h$. It would then make it necessary for the ECF to have crossed 0 \emph{between} the two clusters.  Before stating the theorem, we need a couple of Lemmas first. 
\begin{lemma}
\label{lemma1}
For fixed $n\geq1$, suppose $Y_1,\ldots,Y_n$ are i.i.d. random variables satisfying assumptions $A1$ and $A2$ with $E(Y_1^{2+\epsilon})<\infty$ for $\epsilon >0$. Then as $n \rightarrow \infty$,
\begin{equation*}
\frac{1}{\sqrt{n}}E\left(Y_{(n)}-Y_{(1)}\right)\rightarrow 0.
\end{equation*}
\end{lemma}

\begin{lemma}\label{lemma2}
Suppose $Y'_1,\cdots,Y'_k, Y_1,\cdots,Y_n$ are i.i.d continuous random variables satisfying assumptions $A1$ and $A2$ with support over $[\alpha,\beta]$ with $\infty\leq \alpha<\beta\leq \infty$; here $k$ is fixed and does not change with $n$.  Then, as $n \rightarrow \infty$,
\[
P(Y_{(n)}>Y'_{(k)})\rightarrow 1.
\]
\end{lemma}
This now leads us to the theorem which proves why the ECF captures the clusters. 

\begin{theorem}\label{crossover}
Let $k^*$ denote the total random number of observations in the first cluster and suppose that $2 \leq k^* \leq n-1$. Then, 
\begin{enumerate}
\item for $k^*+1 \leq k \leq n-1$, with probability tending to one, 
\begin{equation*}
G_n\left(\frac{k}{n}\right) \leq 0;
\end{equation*}
\item for $k=k^*-1$, with probability tending to one,
\begin{equation*}
G_n\left(\frac{k^*-1}{n}\right) \geq 0.
\end{equation*}
\end{enumerate}
\end{theorem}
Once we have satisfied ourselves of the fact that for large $n$, the ECF crosses over zero between the two clusters, we can proceed to construct the test for jumps based on $p_n$. The central limit theorem for $p_n$  in theorem \ref{normal} provides us with the test for the presence of jumps based on the following observation: suppose we have a sample from the mixture density given in \ref{density}; if there are no clusters amongst the observations, viz., no jumps in the process, then we do not have a mixture density and instead can regard the observations as arising from a normal density with mean $O(n^{-1})$ and variance $O(n^{-1})$. In that case the true split point $p_0$ solving $G(p)=0$ is $1/2$, owing to the symmetry of the normal density.  Our test for jumps or for the presence of clusters should then ascertain, based on the observations, if the empirical split point $p_n$ is `far' away from $1/2$---if this is the case, then the test signifies the presence of clusters or jumps. More generally, we are interested in testing if the observed sample path falls in $\Omega^c$ or $\Omega^j$ defined in (\ref{sets}). 

Our test can formally be stated as follows: Define
\[
    S_n =   \sqrt{n}\left(\frac{\delta_n(p_n)(p_n-0.5)}{\sqrt{\eta_n(p_n)}}\right),
\]
where $\eta_n(p_n)(\delta_n(p_n))^{-2}$ is as defined in Proposition \ref{variance}. To choose which of the complementary sets  $\Omega^c$ and $\Omega^j$ the discretely observed path of $X$ on $[0,1]$ at times $0 \leq \Delta_n \leq 2 \Delta_n \leq \cdots \leq 1$ falls in, we employ the following decision rule:
\[
  \zeta(n,t,\alpha) = \left\{
  \begin{array}{l l}
   \text{choose $\Omega^c$}  & \quad \text{if $|S_n|\leq z_{\alpha/2}$};\\
    \text{choose $\Omega^j$}   &\quad \text{otherwise},\\
  \end{array} \right.
\]
where $z_{\alpha/2}$ is the $\alpha/2$ standard normal percentile. Our test statistic $S_n$ is
free of $\mu$ and $\sigma$ and we therefore are not required to estimate them. In view of this, it becomes evident that our test for jumps based on $S_n$ is asymptotically of level $\alpha$. It is pertinent to note that $P(\lceil np_n \rceil=0)=\left(1-\frac{\lambda}{n}\right)^n \approx e^{-\lambda}>0$ for large $n$ and free from $n$. That is, with probability $e^{-\lambda}$, even for large $n$, we would, in our sample, not have any observations from the second component of the mixture distribution. The implication of this is that no matter how large our sample is, any test constructed, can never attain power equal to $1$. This is so, since we are observing the process $X_t$ for $t \in [0,1]$ it might be the case that $X_t$ might not have had any jumps in $[0,1]$ despite comprising of a jump component.
\section{SIMULATIONS AND EMPIRICAL FINDINGS}\label{data}
In this section, we investigate the effectiveness of our test via simulations and its performance on two datasets pertaining to S\&P 500 Index. Comparison is also made with the test proposed by \cite{ait2}, hereafter referred to as ST test. All computations and results presented in this section have been performed using the estimator $\frac{\eta_n(p_n)}{\delta_n^2(p_n)}$ for the asymptotic variance of $p_n$.
\subsection{Simulations}
\subsubsection{Accuracy under $H_0$}
\begin{figure}[!h]
\centering
\begin{subfigure}[ht]{0.45\textwidth}
\centering
\includegraphics[width=\textwidth]{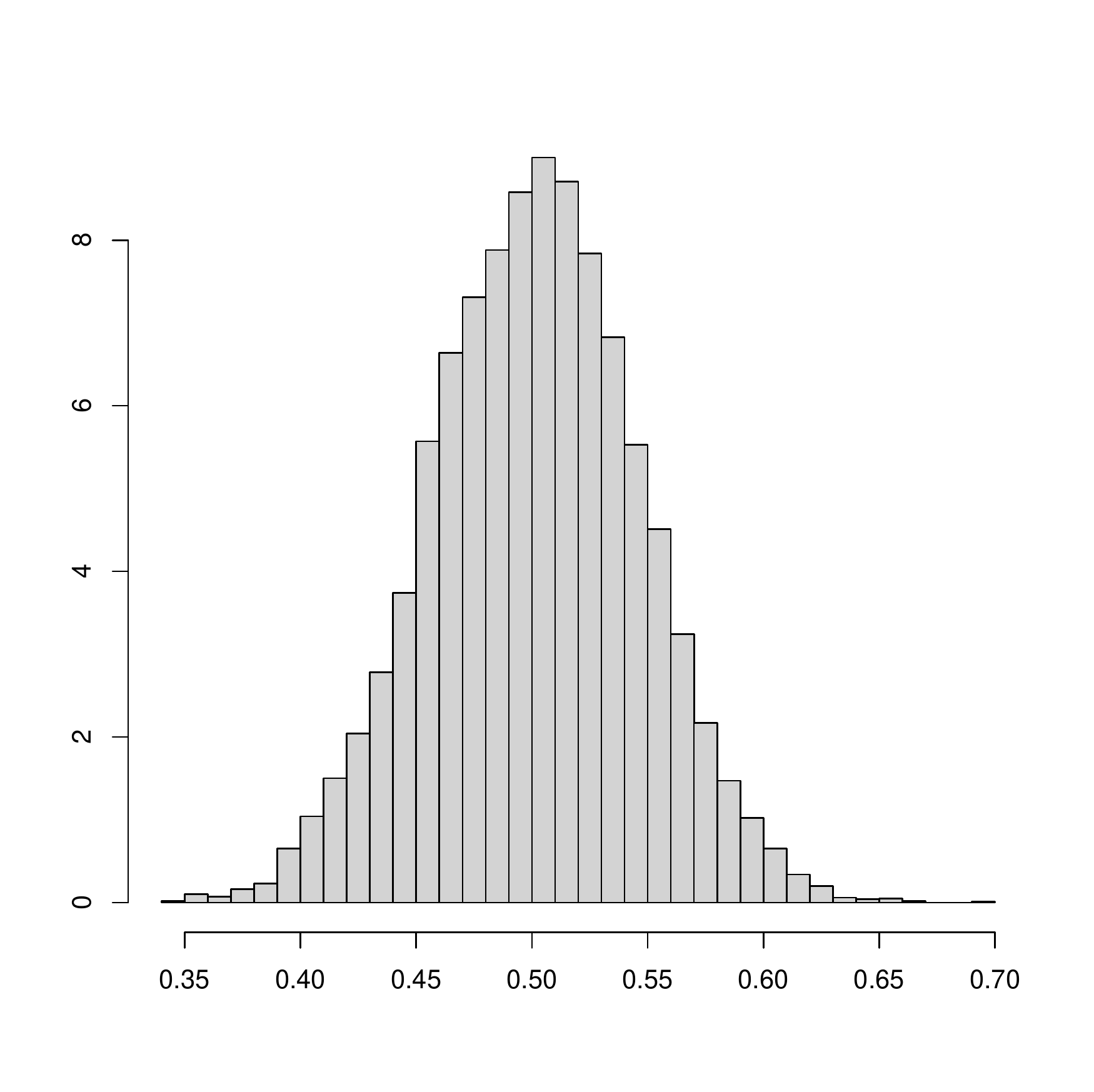}
\end{subfigure}
\centering
\begin{subfigure}[ht]{0.45\textwidth}
\centering
\includegraphics[width=\textwidth]{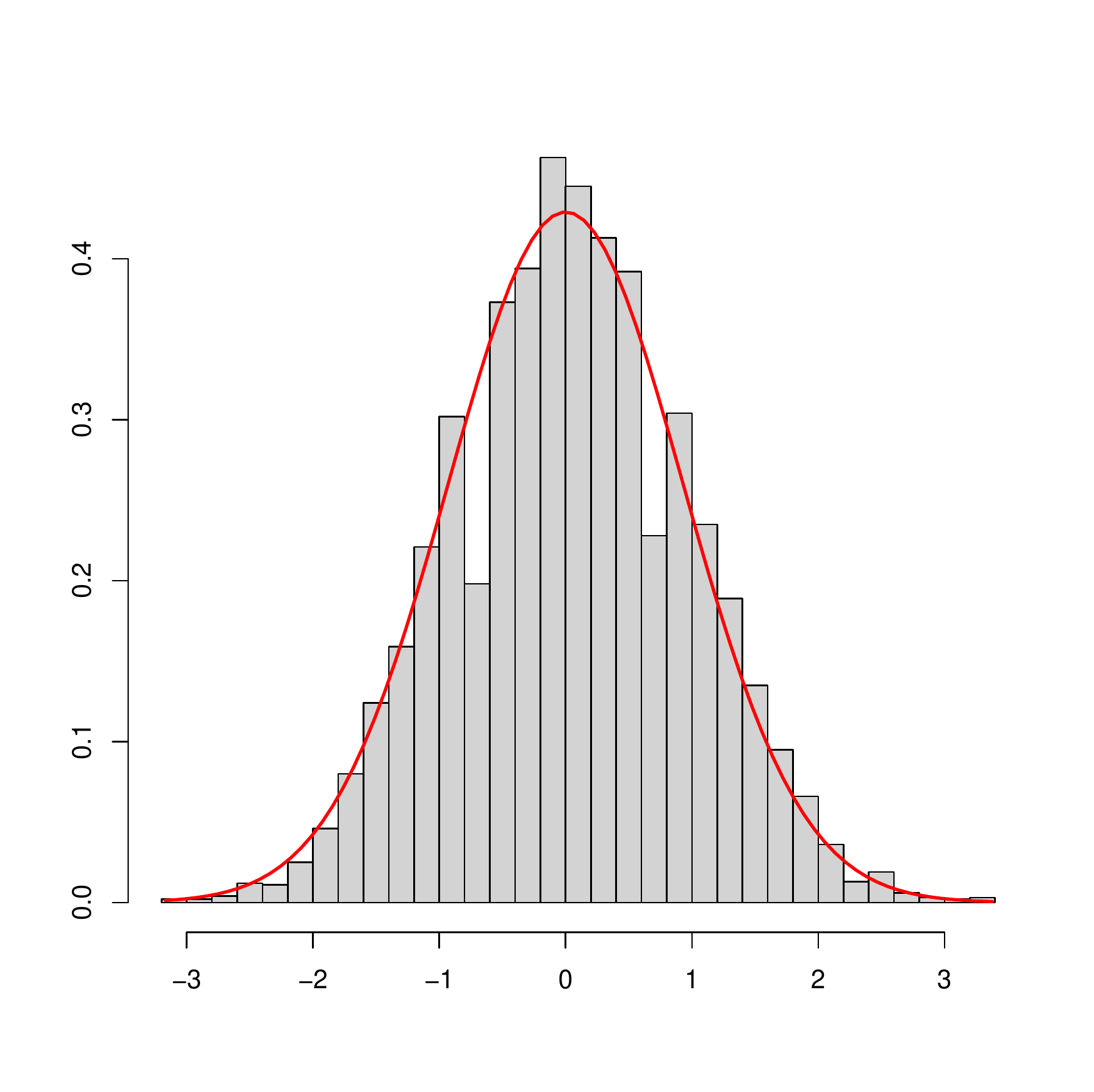}
\end{subfigure}
\caption{\small{Monte-Carlo and theoretical asymptotic normal distribution of the test statistic $p_n$ not standardized (left) and scaled and suitably normalized (right) under the null hypothesis of no jumps. The solid curve on the right represents the theoretical $N(0,1)$ distribution. Here, in order to demonstrate the invariance of $p_n$ to drift and the volatility parameters, we have used  $\mu=4$ and $\sigma^2=1.5$. }}
\label{fig3}
\end{figure}

We first use the theoretical results of $p_n$ in order to set up a simulation of the rejection rate of our test under the null hypothesis of no jumps or in choosing the set $\Omega^c$.  We simulate increments at sampling rate $\Delta_n=1/n$ from a model $dX_t=\mu dt+\sigma d B_t$ containing just the Brownian motion with constant drift $\mu=0$ and spot volatility $\sigma^2=1$; this is done since our test and the ST test are independent of the parameters and remain unaffected by their choice. The test statistic from the ST test in \cite{ait2} given by equation (12) in their paper is computed with $k=2$, $p=4$. The asymptotic level of our test is verified using the CLT for $p_n$ and is compared to the performance of the ST test; it is found that our test requires fewer number of observations, as compared to the ST test, to attain level $\alpha=0.05$. This should not be surprising since the ST test is applicable under a very general setup for a large class of semimartingales. Figure \ref{fig3} depicts the empirical distribution of the non-standardized and standardized test statistic $p_n$. The result from Theorem \ref{normal} appears to be verified by Monte-Carlo simulations. 
\begin{table*}[h]
\begin{center}
\vspace{2mm}
\begin{tabular}{|c|c c |c c| }
\hline
&\multicolumn{2}{c} {\bf{\small{Mean value of $p_n$}}}&\multicolumn{2}{|c|} {\bf{\small{Rejection rate in simulations}}}\\
\hline
$n$&\bf{\small{Asymptotic}} & \bf{\small{Simulations}}&\bf{\small{Our test}}&\bf{\small{ST}}\\
\hline
500& 0.5& 0.484&$0.043$& $0.1037$ \\
1000& 0.5& 0.491&$0.046$& $0.0776$ \\ 
5000& 0.5& 0.503&$0.0492$& $0.0452$ \\ 
10000& 0.5& 0.499&$0.0497$& $0.0418$ \\ 
25000& 0.5& 0.500&$0.0482$& $0.0465$ \\ 
50000& 0.5& 0.501&$0.0501$& $0.0505$ \\ \hline
\end{tabular}
\end{center}
\caption{\small{This table reports the accuracy of $p_n$ as an estimator of $p_0$ for large samples and provides a comparison of level of our test and ST under the null hypothesis of no jumps with $\alpha=0.05$. For the ST test, $p=4$, $k=2$ and $\Delta_n=\frac{1}{n}$ were chosen in order to compute their test statistic. We perform 10000 simulations with $\mu=0$ and $\sigma=1$ since both tests do not depend on them. }}
\end{table*}
\subsubsection{Power against specific alternatives}
In order to examine the power of our test against particular alternatives, we consider two models which differ in their jump components. We consider the model 
\begin{equation}\label{sim_model}
dX_t=\mu dt+\sigma dW_t+J_{t-}dP_t(\lambda) 
\end{equation}
where $P_t(\lambda)$ is a Compound Poisson process with jump sizes $J_t$ given by Double Exponential jumps with location $4$ and scale $1$, proposed in \cite{kou} and \cite{kou2}; we also consider the same model with $J_t$ corresponding to normal jumps with means $1.5$, $10$ and variances $2$, $1$ respectively. Finally, we consider a model with a the jump component driven by a Bernoulli process with success probability $\lambda=0.2$. This model was proposed in \cite{bernoulli} and stands in direct comparison to the Merton model with constant jump sizes and the resulting density in (\ref{density}). Results pertaining to the simulations are in Table 2.

\begin{table*}[h]
\begin{center}
\vspace{2mm}
\begin{tabular}{|c|c c c |c| }
\hline
&\multicolumn{3}{c} {\bf{\small{Compound Poisson process}}}&\multicolumn{1}{|c|} {\bf{\small{Bernoulli process}}}\\
\hline
$n$&\bf{\small{$N(10,2)$ jumps}} & \bf{\small{$DE(4,1)$ jumps}}&\bf{\small{$N(1.5,1)$ jumps}}&\bf{\small{$N(10,2)$ jumps}}\\
\hline
100& 0.036& 0.081& 0.291&$0.039$ \\
1000& 0.048& 0.043& 0.216&$0.047$ \\ 
5000& 0.044& 0.053& 0.187&$0.051$ \\ 
10000& 0.049& 0.047& 0.204&$0.049$ \\ 
25000& 0.052& 0.050& 0.166&$0.048$ \\ \hline
\end{tabular}
\end{center}
\caption{\small{This table reports the results obtained from 10000 simulations pertaining to the rejection rate of our test for two different jump diffusion models given in \ref{sim_model}. The model with the Compound Poisson component is checked with three jump sizes: Normal with mean $10$ and variance $2$, Normal with mean $3.5$ and variance $1$ and Double Exponential (DE) with location $4$ and scale $1$. The jump sizes in the Bernoulli model are Normal with mean 10 and variance 2 with success probability 0.l. Here the Brownian motion drift $\mu=2$ and volatility $\sigma=1$. The tests were conducted at level $\alpha=0.05$ and the corresponding proportion of tests which fail to choose the set $\Omega^j$ of jumps obtained via simulations are reported.}}
\end{table*}
The results in Table 2 further corroborate our assumption of the constant jump size in Section \ref{merton}. Our claim was on the insensitivity of our clustering-based test in practice on the actual nature of the jump sizes as long as they were large jumps; the power of our test against the models with Compound Poisson jump component and Normal or Double exponential jumps is quite good. The Bernoulli jump component is an interesting case in the sense that its compatibility with the density in (\ref{density}) offers a natural setting for the employment of our test. What is clear also is the ordinary performance of our test in the case of small jumps as with the normal jumps with mean 1.5. Our test appears to find it hard to separate out the small jumps since their means is quite close to the Brownian drift $\mu=2$. 

\begin{figure}[!h]
\centering
\begin{subfigure}[ht]{0.45\textwidth}
\centering
\includegraphics[width=\textwidth]{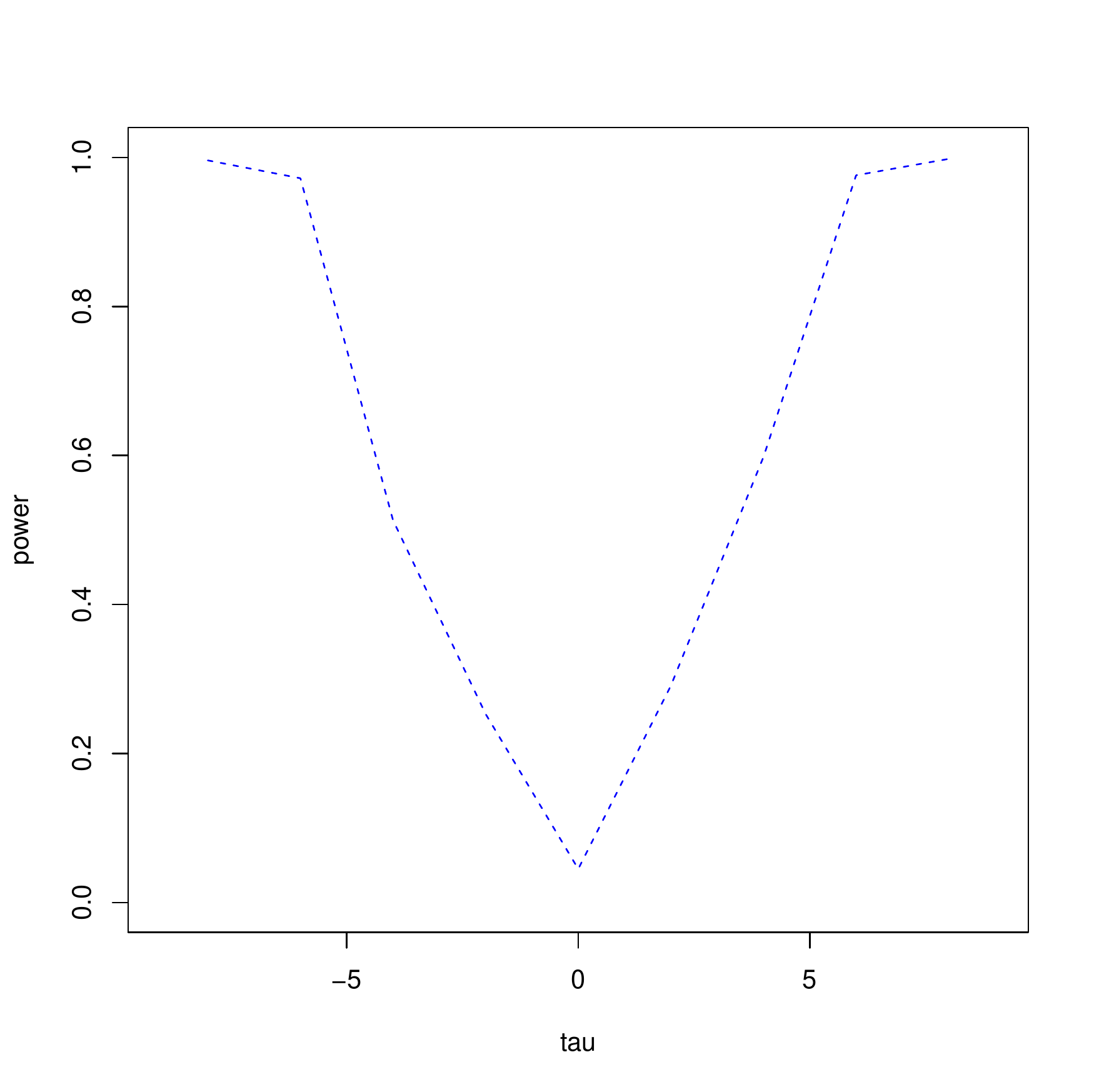}
\end{subfigure}
\centering
\begin{subfigure}[ht]{0.45\textwidth}
\centering
\includegraphics[width=\textwidth]{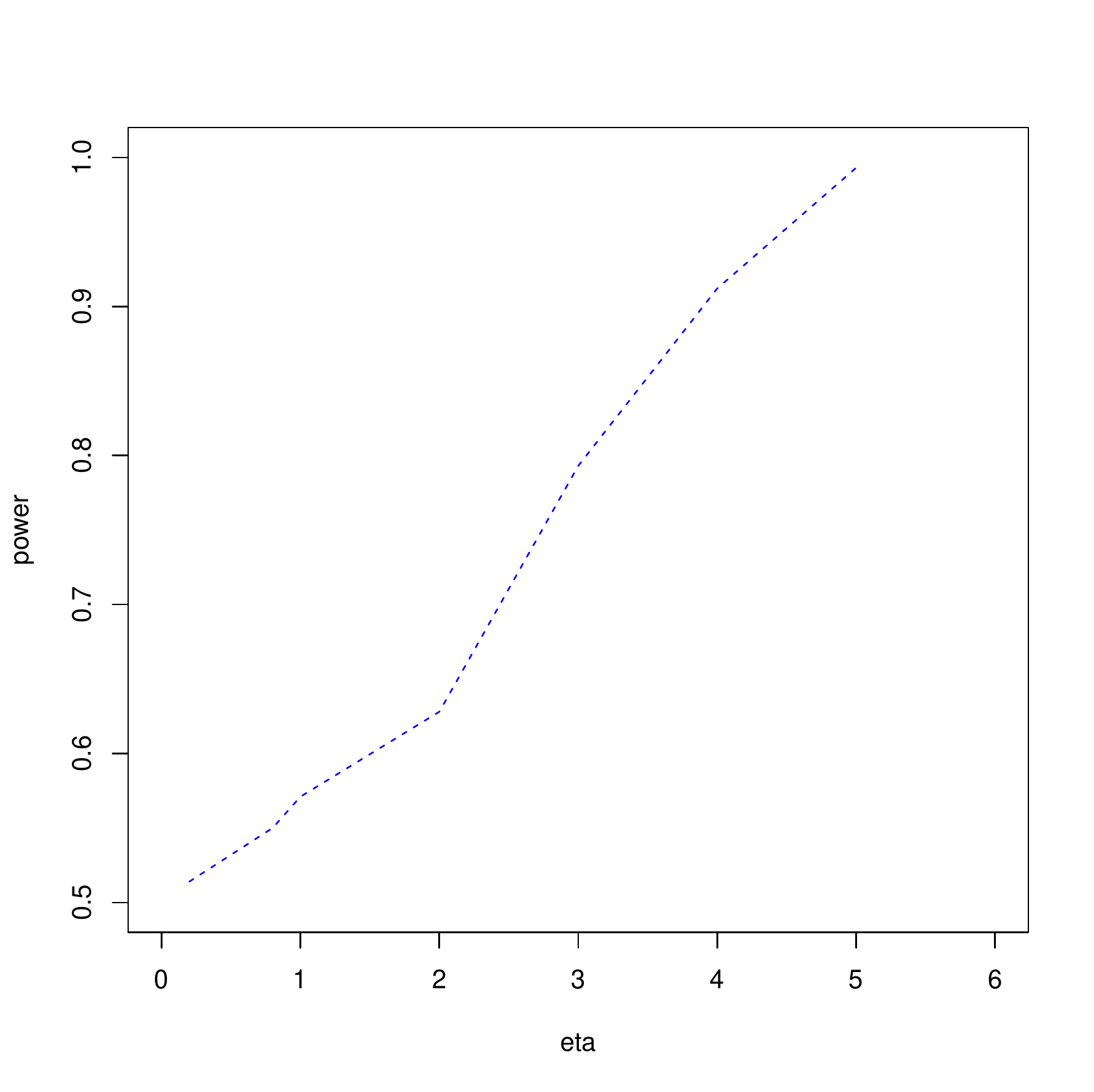}
\end{subfigure}
\caption{\small{Left panel: Power function for our test evaluated at the sequence of alternative models (\ref{sim_model}) by varying the mean of Normal jumps $\tau$; here $\lambda=0.2$, $\eta=1$, $\mu=2$, $\sigma^2=1$. Right panel: Power function for the same model with fixing $\tau=4$ and varying the jump-size variance $\eta$.}}
\label{fig4}
\end{figure}

We now examine the power curves obtained for the model in (\ref{sim_model}) with $P_t(\lambda)$ being a Compound Poisson process with rate $\lambda=0.2$ and $J_t$ are i.i.d. normal random variables with mean $\tau$ and variance $\eta$. Power curves by varying $\tau$ and $\eta$ are provided in Figure \ref{fig4}. The fact that our test is essentially a `story-of-means' is captured in the power functions. Our test, based on the clustering criterion, has good power against models with very large or very small jumps as opposed to models which have small jumps which are comparable to the `jumps' due to the diffusion; this behavior is captured by the fact that when the size of the jumps is similar to the drift, we are required to have high variability in the jump sizes to detect departures from the model with no jumps. 
\subsection{Empirical study}
We conduct the test for two interesting datasets pertaining to daily S\&P 500 Index returns for 500 leading companies which are publicly held on either NYSE or NASDAQ and covers 75\% of U.S equities. Data was obtained from the Federal Reserve Bank of St. Louis Database (\href{http://research.stlouisfed.org/fred2/series/SP500/downloaddata}{research.stlouisfed.org/fred2/series/SP500/downloaddata}) which is publicly available. The first dataset is on index values from Jan 1 1996 to Dec 31 2000 collected daily, wherein there was no discernable jumps in the market; the second dataset is on the index values from Jan 1 2006 to Dec 31 2010, a time window which coincides with the market crash in late 2008. 

\begin{figure}[!h]
\centering
\includegraphics[scale=0.40]{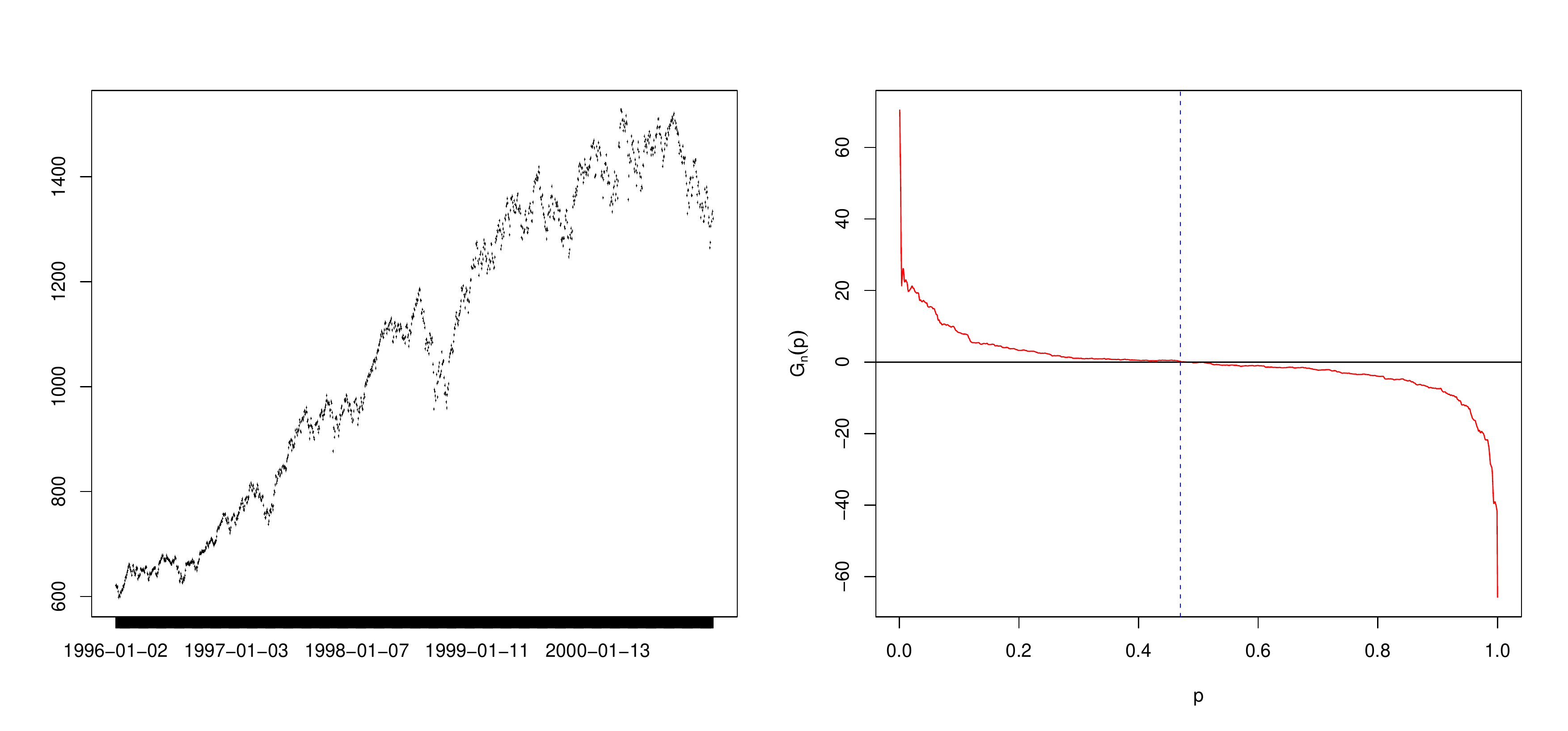}
\caption{\small{On the left is the raw data (S\&P 500 Index) from Jan 1 1996 to Dec 31 2000. On the right is the ECF $G_n$ for the data with $p_n=0.479$ indicated by the dashed blue line. A 95\% asymptotic Confidence Interval for the true split point using the Central Limit Theorem for $p_n$ is $[0.369,0.589]$.}}
\label{fig5}
\end{figure}

For the first data shown in Figure \ref{fig5},  the ECF crosses over zero very close to 0.5 with $p_n=0.479.$ A 95\% confidence interval for the true split point is $[0.369,0.589]$ and since $0.5$ is included in this interval, our test would fail to reject the null hypothesis of no jumps and choose the set $\Omega^c$. This appears reasonable upon viewing the raw index values shown on the left panel of Figure \ref{fig5}. In contrast, for the second dataset the ECF crosses over far away from 0.5 with $p_n=0.237$ with a corresponding 95\% confidence interval being $[0.169,0.305]$ not containing 0.5; this results appears to be consistent with the raw index data which comprises a clear jump downwards corresponding to the market crash in late 2008. In fact, the ST test is in agreement with our test on both the datasets; we omit details in the interests of brevity. 

\begin{figure}[!h]
\centering
\includegraphics[scale=0.40]{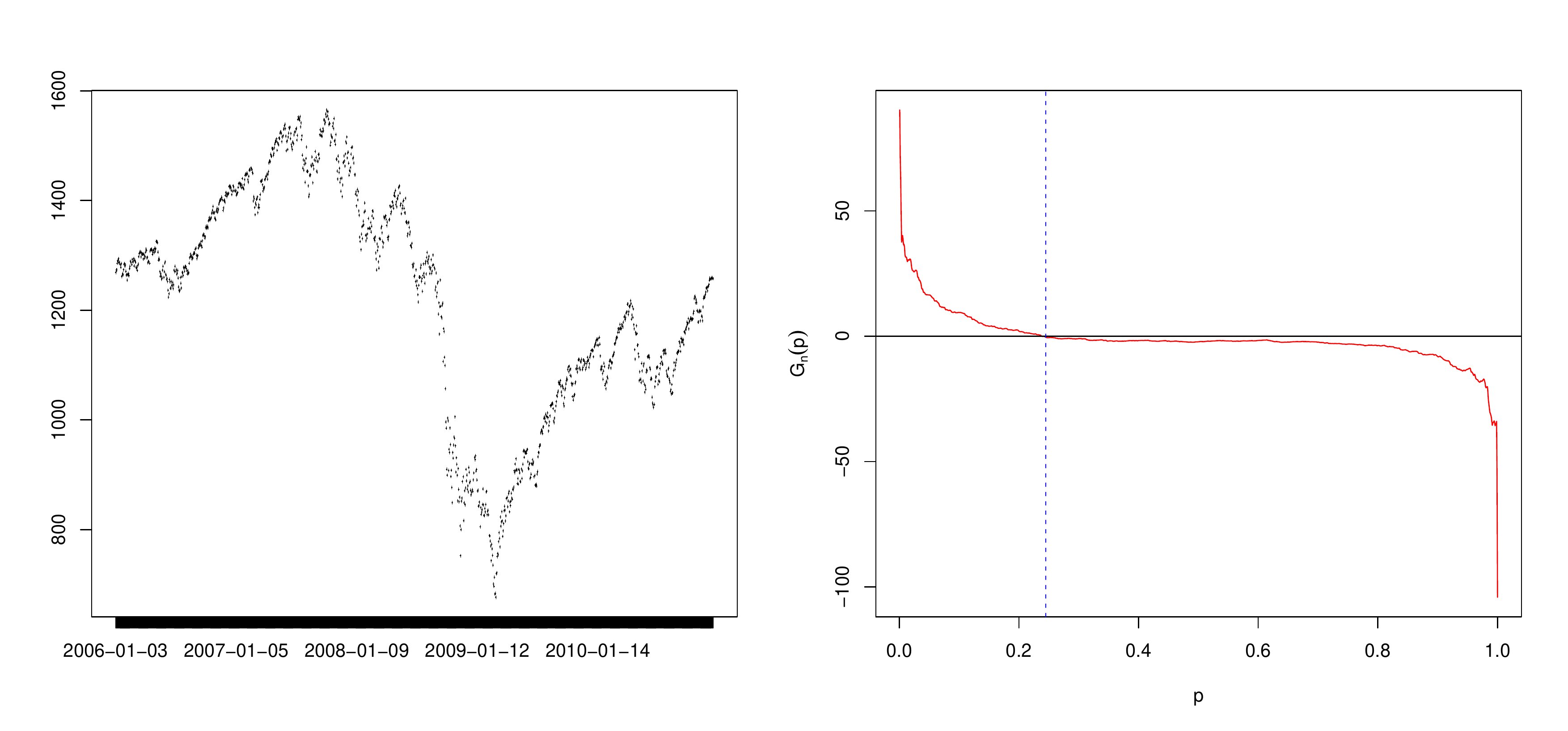}
\caption{\small{On the left is raw data (S\&P 500 Index) from Jan 1 2006 to Dec 31 2010. Following the crash of the market in late 2008 and early 2009, a clear negative jump can be seen in the plot. On the right is the ECF $G_n$ with its corresponding $p_n$ being $0.237$ indicated by the dashed blue line. A 95\% asymptotic Confidence Interval for the true split point result in $[0.169,0.305]$ not containing 0.5. Our test rejects the null hypothesis of no jumps and chooses the set $\Omega_j$.}}
\label{fig6}
\end{figure}
\section{CONCLUDING REMARKS}\label{conclusion}
We have proposed a novel test for the presence of jumps in discretely observed jump diffusion models used in financial applications, based on the simple idea of clustering of observations. The ascendant premise is on the bifurcations of increments of the observed continuous time stochastic process model into those which correspond to the continuous component of the model and those which belong to jumps. While existing methods have concentrated on techniques based on power variations, we have developed our test based on identifying the optimal level of truncation which provides the necessary bifurcation. 

In the case of the popular Merton model for option pricing and its variants, it is shown that the problem of testing for jumps can be reduced to an equivalent one of testing for the presence of clusters in the increments data. Consequently, the asymptotic results from the clustering criterion proposed in \cite{KB2} are used in developing the test. It is to be noted, however, that their results, in existing form, are not directly applicable to the testing problem; one of the contributions in this article is in the modification of the results in \cite{KB2} to suit the requirements of the problem of testing for jumps. 

While the article primarily illustrates the idea of clustering on a special case of the jump diffusion models wherein the drift and instantaneous volatility coefficients are constants, its contribution ought to be viewed as a first step towards investigating the utility in viewing the testing for jumps problem as a testing for clusters problem; indeed, literature is rife with statistical methodologies for handling clustering problems. More importantly, as argued in \cite{kou} and \cite{kou2}, the jump diffusion model is popular amongst practitioners for its analytical tractability and ease of interpretation. In this regard, the results presented here are potentially useful as a quick, perhaps exploratory, check for the presence of jumps; indeed, a simple $k$-means algorithm would suffice. Notwithstanding the absence of generality in the assumed model, it can be noted that the power of the proposed test against various alternatives models which are jump diffusions is promising. The next step would be to consider a more general class of semimartingale models---not unlike the one is \cite{ait2}---and examine the utility of the clustering framework under such a setup. Much work remains to be done in this direction. 
\section{APPENDIX} \label{proofs}
\subsection{Proof of Proposition \ref{variance}:}
Let us first look at the estimator $\delta_n(p_n)$ which is used to estimate $G'(p_0)$, given by
\[
G^{'}(p_0)=\frac{1}{p_0}\left[Q(p_0)-\frac{EW_1\mathbb{I}_{W_1 \leq Q(p_0)}}{p_0}\right]
-\frac{1}{1-p_0}\left[Q(p)-\frac{EW_1\mathbb{I}_{W_1>Q(p_0)}}{1-p_0}\right]-2Q'(p_0).
\]
We will describe in detail how $\delta_n(p_n)$ is consistent for $G'(p_0)$; the argument to show the consistency of $\eta_n(p_n)$ for $Var(\theta_{p_0})$ follows along similar lines once it is noted (following some cumbersome algebra) that $Var(\theta_{p_0})$ is described completely by the following terms: $p_0$, $Q(p_0)$, $f(Q(p_0))$, $\frac{1}{p_0}EW_1\mathbb{I}_{W_1 \leq Q(p_0)}$, $\frac{1}{1-p_0}EW_1\mathbb{I}_{W_1 > Q(p_0)}$,
\begin{equation*}
\frac{1}{p_0}E[W_1^2{\mathbb I}_{W_1<Q(p_0)}], \mbox{ and } \frac{1}{1-p_0}E[W_1^2{\mathbb I}_{W_1\geq Q(p_0)}].
\end{equation*}
We will hence omit the relevant details in that setting. 

Since $p_n$ converges in probability to $p_0$ and using the asymptotic normality of the heavily trimmed sums proved by \cite{stigler}, we have
\[
 T_{nl} \overset{P} \to \frac{EW_1\mathbb{I}_{W_1 \leq Q(p_0)}}{p_0} \quad \text{and} \quad  T_{nu} \overset{P} \to \frac{EW_1\mathbb{I}_{W_1 > Q(p_0)}}{1-p_0}.
\]
Note that 
$$Q^{'}(p_0)=\frac{1}{f(Q(p_0))}$$
is the only other quantity in the expression for $G^{'}(p_0)$ which requires some care with respect to consistent estimation. By assumption $A3$, we have, by definition, 
\begin{align*}
Q^{'}(p_0)=\lim_ {h \to 0}\frac{Q(p_0+h)-Q(p_0)}{h}.\\
\end{align*}

 For $0<p<1$, if $U_{(i)}, 1 \leq i \leq n$ are order statistics corresponding to $n$ $U[0,1]$ random variables, then
\begin{align*}
W_{(\lceil np \rceil +1)}-W_{(\lceil np \rceil)}&\overset{d} = Q(U_{(\lceil np \rceil +1)})-Q(U_{(\lceil np \rceil)})\\
	&\overset{d}=Q(U_{(\lceil np \rceil)}+\Delta_{pn})-Q(U_{(\lceil np \rceil)}),
\end{align*}
where $\Delta_{pn}$ is the uniform spacing $(U_{(\lceil np \rceil +1)}-U_{(\lceil np \rceil)})$. Since we are dealing with central order statistics, i.e. $\frac{\lceil np \rceil}{n} \to p \in (0,1)$, it clear that
\[
\frac{W_{(\lceil np +1\rceil)}-W_{(\lceil np \rceil)}}{\Delta_{pn}} \overset {P}\to \frac{1}{f(Q(p_0))},
\]
as $n \to \infty$, whenever the density $f$ is finite and continuous at the quantile $Q(p_0)$; assumption $A3$ guarantees the fulfillment of these sufficient conditions. In order to obtain $\Delta_{pn}$, note that , for a fixed $0<p<1$, if $F_n$ is the empirical distribution function corresponding to $F$ based on $W_i$, 
\begin{align*}
\Delta_{pn}&=U_{(\lceil np \rceil +1)}-U_{(\lceil np \rceil)}\\
		&\overset{d}=F(W_{(\lceil np \rceil +1)})-F(W_{(\lceil np \rceil)})\\
		&=F_n(W_{(\lceil np \rceil +1)})-F_n(W_{(\lceil np \rceil)})+o_p(1)\\
		&=\frac{1}{n}+o_p(1).
\end{align*}
Furthermore, since $p_n$ is consistent for $p_0$, $W_{(\lceil np_n\rceil)}$ is consistent for $Q(p_0)$ (see p. 308 of \cite{van} for instance). Therefore, a natural estimate of of the inverse of the density at the quantile $Q(p_0)$ would be
\begin{equation}\label{quantile_density}
\hat{Q}^{'}(p_n)=\frac{W_{(\lceil np_n \rceil +1)}-W_{(\lceil np_n \rceil)}}{1/n}.
\end{equation}
Based on the preceding discussion, we can claim that $\hat{Q}^{'}(p_n) \overset{P} \to Q^{'}(p_0)$. Combining these, with a repeated continuous mapping argument, it is easy to note that $\delta^2_n(p_n)$ converges in probability to $G^{'}(p_0)^2$. Putting together these individual pieces along with a continuous mapping argument with the function, it is easy to see note that $\delta_n^2(p_n)$ converges in probability to $G^{'}(p_0)^2$.
\subsection{Proof of Lemma \ref{supGn}:}
Pick any $0<b<1$. Now
\begin{align*}
\sup_{p>b}G_n(p)&=\sup_{p>b}\bigg[\frac{1}{\round}\sum_{j=1}^{\round} W_{(j)}-W_{(\round)}+
\frac{1}{\roundp}\sum_{j=\round +1}^{n} W_{(j)}-W_{(\round +1)}\bigg]\\
& \leq \frac{1}{n}\sum_{j=1}^n W_{(j)}-W_{(\lceil nb \rceil)}+\frac{1}{\lceil n(1-b) \rceil}\sum_{j=\lceil nb \rceil +1}^{n} W_{(j)} -W_{(\lceil nb\rceil )}.
\end{align*}
By assumption $A2$ and the Law of Large Numbers for i.i.d. random variables the first term converges to $0$ in probability. 
Since $W_{(\lceil nb \rceil)}$ is consistent for $Q(b)$, as $b \to 1$, $-W_{(\lceil nb \rceil)} \overset{P}\to -\infty$. As a consequence, the proof of the Lemma would be complete if we can show that 
$$\frac{1}{\lceil n(1-b) \rceil}\sum_{j=\lceil nb\rceil +1}^{n} [W_{(j)}-W_{(\lceil nb\rceil )}]$$
 is bounded in probability for $b$ arbitrarily close to 1. 

Let $\lceil nb \rceil \leq k \leq n$. Now, observe that
\begin{align*}
M_k &=\frac{1}{n-k}\sum_{j=k+1}^n \left[W_{(j)}-W_{(k+1)}\right]\\
&=\frac{1}{n-k}\bigg[(W_{(k+2)}-W_{(k+1)})+\cdots+(W_{(n)}-W_{(k+1)})\bigg]\\
&=\frac{1}{n-k}\bigg[\xi_{k+2}+(\xi_{k+3}+\xi_{k+2})+\cdots+(\xi_n+\xi_{n-1}+\cdots+\xi_{k+2})\bigg]\\
&=\frac{1}{n-k}\bigg[(n-k-1)\xi_{k+2}+(n-k)\xi_{k+3}+\cdots+\xi_n \bigg],
\end{align*}
where $\xi_k=W_{(k)}-W_{(k-1)}$ are the spacings. By the Cauchy-Schwartz inequality, 
\begin{align*}
M_k &= \frac{n-k-1}{n-k}\xi_{k+2}+\frac{n-k-2}{n-k}\xi_{k+2}+\cdots+\xi_n\\
	&\leq \bigg[\left(\frac{n-k-1}{n-k}\right)^2+\cdots+1\bigg]^{1/2}\bigg[\xi^2_{k+2}+\cdots+\xi^2_{n}\bigg]^{1/2}\\
	&\leq (n-k)^{1/2}\big[\xi^2_{k+2}+\cdots+\xi^2_{n}\big]^{1/2}.
\end{align*}
and hence
\[
	\sup_{\lceil nb \rceil \leq k \leq n}M_k \leq \big[\lceil 1 -b\rceil n\big]^{1/2}\big[\xi^2_1+\cdots+\xi^2_{n}\big]^{1/2}.
\]
It is now required that the sum of the squares of the spacings be bounded in probability be of order $n$, in which the expression to the right of the preceding inequality would be of order $(1-b)$; choosing a $b$ close to 1 then completes the proof. This, however, is readily available from Theorem 3 in \cite{PH}, sufficient conditions for which are easily satisfied by the normal distribution. This concludes the proof. 
\subsection{Proof of Lemma \ref{lemma1}:}
\begin{proof}
Denote $\gamma(n,p)=\left\{\displaystyle \int_{0}^1[u^{n-1}-(1-u)^{n-1}]^{p/(p-1)}du\right\}^{(p-1)/p}$ , $R_n=Y_{(n)}-Y_{(1)}$ and $c_{p,n}^p=E|\frac{1}{\sqrt{n}}Y_j|^p$ with $p>2$,.
From theorem $6$ of \cite{BA}, substituting $c_{p,n}^p$ for $\int_0^1|Q(p)|^pdp$,
\begin{equation*}
E(R_n)=\displaystyle \int_0^1n[u^{n-1}-(1-u)^{n-1}]Q(u)du.
\end{equation*}
By H\"{o}lder's inequality, the RHS is
\begin{align*}
	&\leq n\left[\displaystyle \int_0^1|u^{n-1}-(1-u)^{n-1}|^{p/(p-1)}du\right]^{(p-1)/p}\left[\displaystyle \int_0^1|Q(u)|^pdu\right]^{1/p}\\
	&=n\gamma(n,p)c_{p,n}\\
	& = o(1) ,
\end{align*}
for $p>2$, since $\gamma(n,p)=O(n^{-(p-1)/p})$ and $c_{p,n}=O(n^{-1/2})$.
\end{proof}
\subsection{Proof of Lemma \ref{lemma2}:}
\begin{proof}
Fix $\epsilon>0$. Choose $c<\beta$ such that
\[
P(Y'_{(k)}<c)=1-\frac{\epsilon}{2}.
\]
Since $P(Y_i>c)>0$, we have that
\[
	\displaystyle \sum_{i=1}^{\infty}P(Y_i>c)=\infty,
\]
and hence by Borel-Cantelli's lemma,
\[
	P(Y_i>c\textrm{ i.o })=1.
\]
Therefore, there exists $N$ such that
\[
	P(Y_{(n)}>c)=1-\frac{\epsilon}{2} \qquad \textrm{for all } n >> N.
\]
For all $n>>N$, we have
\begin{align*}
	P(Y_{(n)}>Y'_{(k)})&\geq P(Y_{(n)}>c \textrm{ and } Y'_{(k)}<c )\\
&= 1-\frac{\epsilon}{2}+1-\frac{\epsilon}{2}-P(Y_{(n)}>c \textrm{ or }Y'_{(k)}<c )\\
& \geq  1-\frac{\epsilon}{2}+1-\frac{\epsilon}{2} -1 =1-\epsilon.
\end{align*}
\end{proof}
\subsection{Proof of Theorem \ref{crossover}:}
We start first with $1$ which says that $G_n$ is non-positive after it crosses the first cluster comprised of $k^*$ points. For $k^*+1 \leq k \leq n-1$,
\begin{align*}
P\left[G_n (k/n) \leq 0 \right]&=P\left[\frac{\sum_{i=1}^{k^*}W_{(i)}}{k}+\frac{\sum_{i=k^*+1}^{k}W_{(i)}}{k} -W_{(k)}
+\frac{\sum_{i=k+1}^{n}W_{(i)}}{n-k} -W_{(k+1)}\leq 0\right]\\
&=P\left[\frac{k^*}{k}\frac{\sum_{i=1}^{k^*}W_{(i)}}{k^*}+\frac{k-k^*}{k}\frac{\sum_{i=k^*+1}^{k}W_{(i)}}{k} -W_{(k)}
+\frac{\sum_{i=k+1}^{n}W_{(i)}}{n-k} -W_{(k+1)}\leq 0\right]\\
&\leq P\left[\frac{k^*}{k}W_{(k^*)}+\frac{k-k^*}{k}W_{(n)} -W_{(k^*+1)}
+W_{(n)} -W_{(k^*+1)}\leq 0\right]\\
&=P\left[\frac{k^*}{k}\left(W_{(k^*)}-W_{(n)}\right)+2\left(W_{(n)}-W_{(k^*+1)}\right)\leq 0\right]
\end{align*}
Note that the Poisson jump component ensures that the number of observations in the second cluster, $(n-k^*)$ is finite a.s; it is hence the case that $k^*/k\to 1$ as $n \to \infty$. If we can now show that events $\{|W_{(k^*)}-W_{(1)}|>\epsilon\}$ and $\{|W_{(n)}-W_{(k^*+1)}|>\epsilon\}$ tend to zero for an arbitrary $\epsilon$, then the probability of their intersection would go to zero. This would then imply a clear separation of $-h$ between the two clusters and prove that $G_n(k/n) \leq 0$ with high probability for $k^*+1 \leq k \leq n-1$. The probability of both the events tend to 0 owing to Lemmas \ref{lemma1} and \ref{lemma2} respectively.  This concludes the proof of part 1. We now turn our attention to part 2. For $k=k^*-1$, we have
\begin{align*}
P\left[G_n (k/n) \geq 0 \right]&=P\left[\frac{\sum_{i=1}^{k^*-1}W_{(i)}}{k^*-1}-W_{(k^*-1)}+\frac{\sum_{i=k^*}^{n}W_{(i)}}{n-k^*+1} -W_{(k^*)}\geq 0\right]\\
&=P\bigg[\frac{\sum_{i=1}^{k^*-1}W_{(i)}}{k^*-1}-W_{(k^*-1)}-W_{(k^*-1)}+\frac{W_{(k^*)}}{n-k^*+1}\\
&\qquad +\frac{n-k^*}{n-k^*+1}
 +\frac{n-k^*}{n-k^*+1}\frac{\sum_{i=k^*+1}^{n}W_{(i)}}{n-k^*} -W_{(k^*)}\geq 0\bigg]\\
&\geq P\left[(W_{k^*}-W_{(1)})+\frac{W_{(k^*)}}{n-k^*+1}+\frac{n-k^*}{n-k^*+1}W_{(k^*+1)}-W_{(k^*)} \geq 0 \right]\\
&=P\left[ -(W_{(k^*)}-W_{(1)})+\frac{n-k^*}{n-k^*+1}(W_{(k^*+1)}-W_{(k^*)})\geq 0\right]
\end{align*}
A similar argument using Lemmas \ref{lemma1} and \ref{lemma2} as in part 1 concludes the proof. 
\bibliography{ref}
\bibliographystyle{plainnat}
\end{document}